\documentstyle[amssymb,12pt]{amsart}
\title{Sums of Darboux  and Continuous Functions}

\newcommand{\stopproof}{\hfill \nobreak\medskip $\blacksquare$ \\
\hspace*{\fill}}
\newcommand{\CH}{$2^{\aleph_0} = \aleph_1$}

\newcommand{\AND}{\mbox{ \rm and }}

\newcommand{\forces}[2]{\Vdash_{#1} \mbox{``} #2 \mbox{''}}

\newcommand{\proof}{{\bf Proof:} \ }

\newcommand{\RR}{{\Bbb R}}

\newcommand{\SS}{{\Bbb S}}

\newcommand{\lcard}{\, \mid \!}
\newcommand{\rcard}{\! \mid \,}

\newcommand{\card}[1]{\lcard #1 \rcard}
\newtheorem{theor}{Theorem}[section]

\newtheorem{defin}{Definition}[section]
\newtheorem{corol}{Corollary}[section]
\newtheorem{lemma}{Lemma}[section]
\newtheorem{quest}{Question}[section]

\author{Juris Stepr\={a}ns}
\address{Department of Mathematics\\ York University\\ Toronto,
Canada\ \ \ \ \ M3J 1P3 \\{\sc steprans\@nexus.yorku.ca}}

\begin{document}
 \bibliographystyle{amsplain}
 \thanks{This research was partially supported by NSERC}
 \maketitle
 \begin{abstract}
 It is shown that that for every Darboux function $F$ there is a non-constant
 continuous function  $f$  such that
 $F+f$ is still Darboux.
 It is shown to be consistent --- the model used is iterated Sacks forcing ---
 that for every Darboux function $F$
 there is a nowhere constant continuous function                     $f$ such that
 $F+f$ is still Darboux. This answers questions raised in  \cite{be.na.unbd}
 where it  is  shown  that in various models of set
 theory there are universally bad Darboux functions, Darboux functions whose sum with
 any nowhere constant, continuous function fails to be Darboux.
 \end{abstract}

 \section{Introduction}
 A function which maps any connected set to a connected set is known as
 a {\em Darboux function}. This paper will be concerned with functions
 from $\RR$ to $\RR$ and, in this context, Darboux simply means that
 the image of any interval is an interval. While there are various
 results establishing similarities between continuous functions and
 Darboux functions of first Baire class, the fact that it is possible
 to construct Darboux functions by transfinite induction allows all
 sorts of pathologies to exist. For example, transfinite induction can
 be used to construct a Darboux function $F$ such that the function
 $F(x) + x$ is not Darboux
 \cite{rada.darb}.     In \cite{komj.darb} it is shown that if $\cal G$ is a family of
 functions such that $\card{{\cal G}} ^+ < 2^{\aleph_0}$ then there is a Darboux function $F$
 such that $F+g$ is not Darboux for all $g\in \cal G$.
  This result is extended in \cite{be.na.unbd} where it is established, assuming
 certain set theoretic hypotheses, that there exists a {\em universally  bad}
 Darboux function $f\colon \RR\to\RR$ which means that, for every nowhere
 constant continuous $g\colon\bold \RR\to\RR$, $f+g$  does
 not have the Darboux property.  In unpublished work W. Weiss has shown that a
 universally bad Darboux function can be constructed  assuming only the existence of a
 $2^{\aleph_0}$ additive ideal $\cal I$ on $\cal B$, the Borel subsets fo $\RR$, such that
 the Boolean algebra ${\cal B}/{\cal I}$ has the $2^{\aleph_0}$ chain condition; in other
 words, there do not exists $2^{\aleph_0}$  elements of      $\cal B$ whose pairwise
 intersections belong to $\cal I$.

 In this paper it will be shown that some form of  set theoretic hypothesis
  is necessary for such a result  because there is a model of set theory where for every
 Darboux function $F$ there is a nowhere constant continuous function $F$
 such that $F+f$ is also Darboux. The
 significance of the adjective ``nowhere constant'' in this statement
 requires some comment because it might seem a minor point. An indication that this is not so is given
  by the fact that, in
 spite of having shown that there is a Darboux function $f$ such     for every nowhere
 constant continuous $g\colon\bold \RR\to\RR$, $f+g$  does
 have the Darboux property, the authors of          \cite{be.na.unbd}
 pose the following question at the end of their paper.
 \begin{quest} Does there exist a   Darboux function $F\colon
 \RR \to \RR$ such that $F+g$   does not have the Darboux property whenever
 $g$ is continuous
 but  not constant?
	  \end{quest}

 Section 2 provides a negative answer to this problem. Section 3
 contains some technical material on Sacks forcing             and Section 4 makes use of
 this material in proving the main consistency result.
						       The final section contains some open
  questions.
  \section{Sums with
 Non-constant Functions}
 The following lemma is easy but the  proof is included anyway. It is
 essentially the Sen-Massera condition which, in the case of Baire
 class 1 functions, is equivalent to being Darboux \cite{bruc.deri}.
 \begin{lemma}
  If $F$  is Darboux and not \label{osc}continuous at $x$ then there is an
  interval
  $(a,b)\neq \emptyset$ such that for each $y\in (a,b)$ there is a sequence $\{x_n\mid n\in
 \omega\}$
  such that $\lim_{n\rightarrow\infty}x_n = x$ and $F(x_n) = y$ for all $n\in \omega$.
  \end{lemma}
  \proof    Because $F$ is not continuous at $x$ there are sequences $\{y^a_n\}_{n\in\omega}$
 and $\{y^b_n\}_{n\in\omega}$ such that $$\lim_{n\to\infty} y^a_n = x = \lim_{n\to\infty}
 y^b_n \AND
 \lim_{n\to\infty} F(y^a_n) = a < b \lim_{n\to\infty} F(y^b_n)$$
 Given $n\in \omega$ and
 $y\in (a,b)$   let $k$ be such that $\mid y^a_k - x\mid < \frac{1}{n}$,
 $\mid y_k^b - x\mid < \frac{1}{n}$ and
 $F(y^b_k) > y  > F(y^a_k)$. Then use the Darboux property of $F$ to find $x_n$ between
 $y^a_k$ and $y^b_k$ such that $F(x_n) = y$.
 \stopproof
 \begin{corol}
 If $F$ is a Darboux function \label{ft1} which is  finite-to-one then $F$ is continuous.
 \end{corol}
 \begin{lemma}
 If $F: \RR \to \RR$ is a Darboux function which is\label{nowcts}  continuous at
 only countably many points then there is
 a non-constant, continuous function $f$ such that $F+f$ is Darboux.
 \end{lemma}
 \proof  To each real $x$ at which $F$ is not continuous,
 use Lemma~\ref{osc} to assign an
 interval $(a_x, b_x)$
 such that for each $y\in (a_x,b_x)$ there is a sequence $\{x_n\mid n\in \omega\}$
  such that $\lim_{n\rightarrow\infty}x_n = x$ and $F(x_n) = y$ for all $n$. For rationals
 $p$ and  $ q$
 let $X(p,q)$ be the set of all $x$ such that $a_x < p < q < b_x $ and note that
 $X(p,q)$ is a closed set.
 Because $F$ is continuous at only countably many  points, it
 can not be the case
 that $X(p,q)$ is nowhere dense for each pair of rationals
 $p$ and $q$. Therefore let $[s,t]$ and $[a,b]$  be
 intervals such that $[a_x,b_x] \supseteq [a,b]$
 for each $x\in [s,t]$ and, furthermore $F(s) = a$ and $F(t) = b$. Observe that
 $F^{-1}\{y\}$ is dense in $[s,t]$ for each $y \in [a,b]$.

 Next, choose a family of open
   intervals $\cal I$ such that
   \begin{itemize}
 \item ${\cal I} = \cup_{n\in\omega}{\cal I}_n$ where ${\cal I}_n = \{(u^n_i,v^n_i) \mid i\in 2^n-1\}$
 \item $v^{n+1}_i < u^n_i < v^n_i < u^{n+1}_{i+1}$ for
 each $n$ and $i$
 \item $\cup\cal I$ is dense in $[s,t]$
 \item $F(u^n_i) = a$ and $F(v^n_i) = b$ for each $n$ and $i$
 \item $\sup F\restriction [v^{n}_{i-1},u^{n+1}_{i}]\cup [v^{n+1}_{i}, u^{n}_{i}]
 < \sup F\restriction [u^{n+1}_{i},v^{n+1}_{i}] + \frac{1}{n+1}$ where it is understood that, in the case $i = 0$,
 $v^n_{-1} = s$ and, in the case $i = n$, $  u^n_n = t$
 \item $\inf F\restriction [v^{n}_{i-1},u^{n+1}_{i}]\cup [v^{n+1}_{i}, u^{n}_{i}]
 > \inf F\restriction [u^{n+1}_{i},v^{n+1}_{i}] -  \frac{1}{n+1}$ where it is understood that, in the case $i = 0$,
 $v^n_{-1} = s$ and, in the case $i = n$, $  u^n_n = t$
       \end{itemize}

 Now let $g$ be a nondecreasing, continuous function such that
 $g\restriction I$ is constant for each $I\in \cal I$ and such that
 $g(s) = 0$ and $g(t) = (b-a)$ and such that $g$ is constant on
 $(-\infty,s]$ and $[t,\infty)$.          The reader who insists on concreteness may verify
 that $$g(x) = \frac{\sup\{\frac{2i+1}{2^n} \mid u^n_i \leq x\}}{b-a}$$
 satisfies these requirements.

 To see that $F+g$ is Darboux suppose that $x < y$ and $w $ lies between
 $F(x) + g(x)$ and $F(y) + g(y)$.
 First of all observe that it may be assumed that
  $s \leq x < y \leq t$. The reason  it may be assumed that $s\leq x$ is that
 if  $x < s$ then either $ w$ lies between
  $F(x) + g(x)$ and $F(s) + g(s)$ or else it lies between $F(s) + g(s)$ and $F(y) + g(y)$.
  In order to eliminate the first case use the fact that $F$ is Darboux and $g$ is
 constant on $[x,s]$ to find
  $z\in [x,s]$ such that $F(z) +g(z) = w$. In the second case it may, of course, be assumed
 that $x=s$. A similar argument can be applied to show that, without loss of generality, $y
 \leq t$.

 First consider the case that there is some $(u,v)\in \cal I$
 such that  $g$ has constant value $c$ on
 $(u,v)$ and $w \in (a+c ,b+c)$ and such that $x \leq u < v \leq y$.
 It is possible to use the fact that $g$ is constant on $[u,v]$, $F$ is Darboux, $F(u) = a$
 and $F(v) = b$ to find $ z\in (u,v)$ such that
   $F(z) + g(z) = w$.

  In the remaining case it follows from the  fact that  $\cup\cal I$ is dense and the
 continuity of $g$ that either $w \geq b+g(r)$ for every $r \in \cup{\cal I}\cap (x,y)$ or
  $w \leq a+g(r)$ for every $r \in \cup{\cal I}\cap (x,y)$.
   Only the first case  will be considered since the other one is dealt with similarly.
 Furthermore, it will be assumed that $F(x) + g(x) < F(y) + g(y)$ since the other case is also similar.
  To begin,
  suppose that $y \in (u,v) \in \cal I$. Since $F(u) = a < b < w - g(y)$ it follows that
 $F(u) < w - g(y) < F(y)$ and so it is possible to appeal to  the Darboux property of $F$
 and the constantness of $g$ on $[u,v]$.

  On the other hand, if $y\notin
  \cup{\cal I}$ then
  there must be some $m\in \omega$ such that
 $F(y) + g(y)  - w > 1/m$. Choose $\delta > 0$ such that $y-x > \delta$ and
 $g(y) - g(r) > \frac{1}{2m}$ if $y-r < \delta$.
 Then there is
 some  $k\geq 2m$ such that $(u^k_i,v^k_i) \in {\cal I}_k$  and
 $ y - \delta  \leq u^k_i < v^k_i < y\leq u^{k-1}_i $ for some $i\in 2^n-1$. It follows
 that $F(y) \leq \sup F\restriction  [v^{k}_{i}, u^{k-1}_{i}]
 < \sup F\restriction [u^{k}_{i},v^{k}_{i}] + 1/k$ and hence
 $\sup \{F(r) + g(r)\mid r\in [u^{k}_{i},v^{k}_{i}]\} > F(y) + g(y) -1/m$.
  Therefore $w < \sup \{F(r) + g(r)\mid r\in [u^{k}_{i},v^{k}_{i}]\} $. Because
 of the case being considered, it follows that $w \geq b+ g(u^k_i) > a+g(u^k_i) = F(u^k_i) + g(u^k_i)$;
 in other words there is some $r\in [u^{k}_{i},v^{k}_{i}] $ such that
   $$F(u^k_i) + g(u^k_i) < w   < F(r) + g(r)$$
 Because $g$ is constant on $[u^{k}_{i},v^{k}_{i}]$ it
 now follows from the Darboux property of $F$ that there is $z\in
 [u^{k}_{i},v^{k}_{i}]$ such that $F(z) +g(z) = w$.
 \stopproof

       \begin{lemma} If  $F$ is a Darboux \label{ctsm}      function which is continuous on an uncountable
       set  then there is a continuous, non-constant function $g$ such that $F+g$ is Darboux
						      \end{lemma}
 \proof
	   Because the set  of points where $F$ is continuous is Borel, it is possible to find a
	   perfect, nowhere dense set $P$ such that $F$ is continuous at each point of $P$.
 Because $P$ is a  perfect, nowhere dense set it follows that $\RR\setminus P =
 \cup\cal I$ where $\cal I$ is a disjoint family
 of open intervals of order type the rationals. Let $g$ be any continuous, non-decreasing
 function which is not constant yet,  $g$ has constant value $g_I$ on
 each interval $I \in \cal I$.

 To see that $F+g$ is Darboux suppose that $x < y$ and that  $F(x) +g(x) < w < F(y) + g(y)$
 --- a similar proof works if $F(x) +g(x) > w > F(y) + g(y)$. If there is
 some interval $I \in \cal I$ such that
 $I\subseteq [x,y]$ and $\sup F\restriction I \geq w - g_I \geq \inf F\restriction I$
 then the Darboux property of $F$ guarantees that there is
 some $z \in I$ such that $F(z) +g(z) = w$.

  If there is no such $I$ then  consider first the case that there are $I$ and $J$ in
 $\cal I$ such that $I\subseteq [x,y]$ and $J\subseteq [x,y]$ and
  $\sup F\restriction I < w - g_I$ and $\inf F\restriction J > w - g_J $ and suppose
 that $\sup I <\inf J$.
  Let $z$ be the infimum of all intervals $J'$ such that $\sup I <\inf J'$ and
 such that
  $w - g_{J'} <  \inf F\restriction J'$.
  First observe that  $z \notin \cup\cal I$ and so   $F$ is continuous at $z$ and hence
  $$w-g(z) \leq \liminf_{r\to z^+}F(r) = \lim_{r\to z^+}F(r) = F(z)$$
  On the other hand, since  $g$ is also  continuous at $z$ it follows from the
  defining property of $z$ that
   $$w-g(z) \geq \liminf_{r\to z^-}F(r) = \lim_{r\to z^-}F(r) = F(z)$$
							and so $F(z) + g(z) = w$.
  Similar
 arguments in the other cases establish that one the following two possiblities holds
 \begin{itemize}
 \item if $I\subseteq [x,y]$ then  $\sup F\restriction I < w - g_I $
 \item if $I\subseteq [x,y]$ then  $\inf F\restriction I > w - g_I $
 \end{itemize}
 Consider the first alternative. If $y\notin \cup\cal I$ then $F$ is continuous at
 $y$
 and so $$\lim_{s\to y^{-}}F(s) + g(s) = F(y) + g(y) \leq w$$ and  this is impossible
 becasue $F(y) + g(y) > w$. If
 $y \in (a,b) \in \cal I$
 and $F(a) + g(a) > w$ then the same argument
 applies because $a  \notin \cup\cal I$ . On the other
 hand, if  $F(a) + g(a) \leq w$ then the Darboux property of $F$ and the fact
 that $g$ is constant on $[a,b]$
 yields $z\in [a,y)$ such that $F(z) + g(z) = w$. The other alternative is dealt with similarly.
 \stopproof

 \begin{theor}
	 If $F$ is a Darboux function then there is a non-constant
 continuous function $g$ such that
 $F+g$ is Darboux.\end{theor}
 \proof Either $F$ is continuous on an uncountable set or it is not. If it is, use
 Lemma~\ref{ctsm} and if it is not then use Lemma~\ref{nowcts}. \stopproof

 \section{Sacks Reals}
 The Sacks partial order of perfect trees will be denoted by $\SS$ and the
 iteration, of length $\xi$, of
 this partial order will be denoted by $\SS_\xi$ --- so $\SS_1 =
 \SS$ and $\SS_0 = \emptyset$.
 For other  notation and definitions concerning Sacks
 reals see \cite{mill.mapp} as well as \cite{ba.la.sack}.  For  any $p\in
 \SS_{\xi}$ define $$p^* = \{\theta : \xi\times\omega\rightarrow 2 \mid
 (\forall F \in [\xi]^{<\aleph_0})(\forall
 m\in\omega)(\theta\restriction F\times m \mbox{ is consistent
 with } p)\}$$ It is easy to see $p^* \subseteq 2^{\xi\times
 \omega}$ is a closed set; but there is no reason to believe that it should     be
 non-empty.  However, if $p$ is  determined ---  see page 580 of
 \cite{mill.mapp} for a definition --- then $p^*$ is a reasonably
 accurate reflection of $p$. In \cite{mill.mapp} a notion very similar to $p^*$ is defined
 and denoted by $E_p$. The only difference is that $E_p\subseteq (2^\omega)^A$ where $A$
 is the domain of $p$.     The projection function from
 $2^{\xi\times\omega}$ to $2^{\gamma\times\omega}$
 will be denoted by $\Pi_{\xi,\gamma}$.

 \begin{lemma}
 If $p \in\SS_\xi$ is $(E,k)$-determined and $p\forces{\SS_{\alpha}}{x \in \RR\setminus V}$
 then for each $E\in [\alpha]^{<\aleph_0}$ and $k\in\omega$ there is
 $q$ such that $(q,k) \leq_E (p,k)$ and \label{name} a function $Z : q^* \rightarrow
 \RR$ such that
 \begin{enumerate}
 \item $q\forces{}{x = Z(G)}$
 \item $Z(x) \neq Z(y)$ unless $\Pi_{\omega_2,1}(x) = \Pi_{\omega_2,1}(y)$
 \end{enumerate}
 \end{lemma}
 \proof This is essentially Lemma~6 on page 580 of \cite{mill.mapp}. The only difference is that
 it is now required that $(q,k) \leq_E (p,k)$ whereas Miller's Lemma~6
 only asserts that $q \leq p$. On
 the other hand, the assertion required here is only that
 $Z(x) \neq Z(y)$ unless $\Pi_{\omega_2,1}(x) = \Pi_{\omega_2,1}(y)$;
 whereas a canonical condition for $x$, in Miller's terminology, actually yields a
 one-to-one function
 $Z$. The way around this is to choose for each $\sigma : E\times k \rightarrow 2$ a condition
 $q_\sigma$ and a one-to-one function $Z_\sigma : ( q_\sigma\restriction \beta(\sigma))^*
 \rightarrow \RR$
 such that  $q_\sigma \forces{}{x = Z_\sigma(\Pi_{\omega_2,\beta(\sigma)} (G))}$.  The
 point
 to notice is that the domain of $Z_\sigma$ depends on $\beta(\sigma)$ and so there may
 not be a single ordinal
 which works  for all $\sigma$. Nevertheless, $\beta(\sigma) \geq 1$ for each $\sigma$ and so
 it is possible to define $Z = \cup_{\sigma} Z_\sigma\circ \Pi_{\omega_2,\beta(\sigma)}$. It follows
 that $Z(x) \neq Z(y)$ unless $\Pi_{\omega_2,1}(x) = \Pi_{\omega_2,1}(y)$.
 \stopproof
 \begin{lemma}
 If $p \in \SS_\xi$ is $(E,k)$-determined and
 $F:p^*\to \RR$ and  $G:p^*\to \RR$ are   continuous functions such that
  $F\restriction q^*\neq G\restriction q^*$ for each $ q\leq p$ then
 there is some
  $q$ such that $(q,k) <_{E} (p,k)$ and the images of $q^*$ under
 $F$ and $H$ are \label{splitimage}disjoint.
 \end{lemma}
 \proof Let $\Sigma$ be the set of all $\sigma : E\times k \rightarrow 2$ which are
 consistent with $p$. For each $\sigma \in \Sigma$
 $F\restriction (p\mid\sigma)^*\neq G\restriction (p\mid\sigma^*$ and so it is possible
 to find some $x_\sigma \in (p\mid\sigma)^*$, $E_\sigma\in [\xi]^{<\aleph_0}$ and $k_\sigma\geq k$ such
 that
 $F(x_\sigma) \neq G(x_\sigma)$ and, moreover,
 the image of $(p\mid x_\sigma\restriction E\sigma)\times k_\sigma)^*$ under $F$
 is disjoint from the image
 under $G$. Let $q' = \cup_{\sigma \in \Sigma}p\mid x_\sigma\restriction
 E\sigma)_\sigma\times k\sigma)$.
 By repeating this operation for each pair
 $\{\sigma,\tau\}\in [\Sigma]^2 $ it is possible to obtain $q$ with the
 desired properties.
 \stopproof

 \section{Darboux Functions and the Sacks Model}

 \begin{lemma}
 If $H: I \rightarrow \RR$ is Darboux then there is a countable set $D$
 such
 \label{countable_witness} that, for any continuous function $F$, if for
 every $a\in D$ and $b
 \in D$ and $t$ such that $$H(a) + F(a) < t < H(b) + F(b)$$ there is some
 $c$ between $a$ and $b$ such that $H(c) + F(c) = t$ then $H+F$ is also
 Darboux.
 \end{lemma}
      \proof Let $D$ be any countable set such that $(H)\restriction D$ is dense
      in the graph of $H$ and suppose that $F$ is continuous.
  If $F(x) + H(x) < t < F(y) + H(y)$ then, because $F$ is continuous,
      there is some $\epsilon > 0$
      such that $F(z) + H(x) < t $ if $\mid z - x \mid < \epsilon$
      and  $F(z) + H(y) > t $ if $\mid z - y \mid < \epsilon$. Because
 $H$ is Darboux
 and $H\restriction D$ is dense in the graph of $H$
  there are $d_x\in D$ and $d_y \in D$, between $x$ and $y$,
      such that   $\mid d_x - x \mid < \epsilon$  and $\mid d_y - y
 \mid < \epsilon$
 and $H(d_x) < t - H(x)$ and $H(d_y) > t _ H(y)$.
      Hence $F(d_x) + H(d_x) < t $   and $F(d_y) + H(d_y) >  t$.
 Hence, if there is some $z$ between $d_x$ and $d_y$ such that $H(z) +
 F(z) = t$ then $z$ also lies between $x$ and $y$.
 \stopproof

 For the rest of this section by a condition in $\SS_\xi$ will be meant a determined
 condition. Real valued functions will be considered to have as their
 domain the unit interval $I$. This is merely a convenience that allows
 the use of the complete metric space of al continuous real valued
 functions on the unit interval using the $\sup$ metric. This space
 will be denoted by ${\cal C}(I,\RR)$ and its metric will be $\rho(f,g)
 = \sup\{\mid f(x) - g(x) \mid : x\in I\}$.

 \begin{theor}
 Let $V$ be a model of \CH and $V[G]$ be obtained by
 adding $\omega_2$ Sacks reals with countable support iteration.  If $H
 : I \rightarrow \RR$ is a Darboux function in $V[G]$ then
 there\label{main} is a second category set of continuous
 functions $f$ such that $H+f$ is also Darboux.
 \end{theor}
 \proof   If the theorem fails then, in $V[G]$,  let $H$ be a Darboux function
 and $X$ be  a comeagre subset of ${\cal C}(I,\RR)$
  which provide a counterexample.  Let $D$ be a countable set, provided by
 Lemma~\ref{countable_witness}, such that for any continuous function
 $F$, if for every $a\in D$ and $b
 \in D$ and $t$ such that $H(a) + F(a) < t < H(b) + F(b)$ there is some
 $c$ between $a$ and $b$ such that $H(c) + F(c) = t$ then $H+F$ is also
 Darboux.    It must be true
  that, for each continuous function $g\in X$
  there is an interval $N(g) = [a,b]$, with endpoints in $D$,
 and a real $T(g)$ between $H(a) + g(a)$ and $H(b) +
 g(b)$ such that there is no $z\in [a,b]$ such that $H(z) + g(z) =
 T(g)$.

 By a closure argument, there must      exist $\alpha \in\omega_2$
 such that \begin{itemize}
 \item $D\in V[G\cap \SS_\alpha]$
 \item  $T(f) \in V[G\cap \SS_\alpha]$ for every $f\in V[G\cap \SS_\alpha]$
 \item if $x$ is in $V[G\cap \SS_\alpha]$ then so is  $H(x)$
 \item  $X = \cap_{n\in\omega}U_n$ where each $U_n$ is a dense
 open set belonging to $V[G\cap \SS_\alpha]$.
 \end{itemize}

 To simplify notation it may be assumed that $V = V[G\cap \SS_\alpha]$.
 In $V$, let $\{d_i \mid i\in \omega\}$ enumerate $D$,
 let $G$ denote the generic function from $\omega_2\times \omega$ to 2
 which is obtained from an $\SS_{\omega_2}$ generic set and, let $p_0\in \SS_{\omega_2}$
 be a determined condition.

 Let
 $\frak M$ be a countable elementary submodel of $(H(\omega_3),\in)$ containing
 the functions  $T$ and $N$ and the name $H$. Let $\{E_n \mid n\in\omega\}$ be an
 increasing sequence of finite sets such that $\cup_{n\in\omega}E_n =
 {\frak M}\cap \omega_2$. (The use of the elementary submodel is only a
 convenience that allows the finites set $E_n$ to be chosen before
 beginning the fusion argument, thereby avoiding some bookkeeping.) Construct, by induction on
 $n\in\omega$, functions $f_n$, as well as conditions $p_n \in
 \SS_{\omega_2}$, reals $\epsilon_n > 0$ and integers $k_n$, all in $\frak M$, such that:
 \begin{itemize}
 \item[IH(0)] $f_n \in {\cal C}(I,\RR)$ and $f_0$ is  chosen arbitrarily
 \item[IH(1)] the neighbourhood of $f_n$ of radius $\epsilon_n$ in ${\cal C}(I,\RR)$
 is contained in $\cap_{i\leq n}U_n$
 \item[IH(2)] $\rho( f_n - f_{n+1})  < \epsilon_n\cdot 2^{-n-1}$
 \item[IH(3)] $ p_n$ is $(E_n,k_n)$ determined
 \end{itemize}
 For
 each $n$, an integer $J_n$ and a sequence ${\cal C}_n = \{c_j^n \mid
 j\leq J_n \}$ such that $c_0^n = 0 < c_1^n < c_2^n \dots < c^n_{J_n} =
 1$ will be chosen so that
 \begin{itemize}
 \item[IH(4)] $d_n \in {\cal C}_n$ and ${\cal C}_i \subseteq {\cal
 C}_n$ if $i\in n$
 \item[IH(5)] if $i\in n$ and $c\in {\cal C}_i$    then $f_n(c) = f_i(c)$
 \end{itemize}
 For each $n$ and each $j\leq J_n$ a continuous function $\Phi_{n,j} : p_n^*\to \RR$
 will be found so that there is a name $z_{n,j}$ such that
 \begin{itemize}
 \item[IH(6)]  $p_n\forces{\SS_{\omega_2}}{H(z_{n,j}) =
 \Phi_{n,j}(G)}$ for each $j\in J_n$
 \end{itemize}
 A function $Z_{n,j} :2^{\omega_2\times\omega}\rightarrow \RR$ will
 also be constructed so that
 \begin{itemize}
 \item[IH(7)] $p_n
 \forces{\SS_{\omega_2}}{z_{n,j} =
 Z_{n,j}(G)}$
 \item[IH(8)] if $Z_{n,j}(x) = Z_{n,j}(y)$ then $\Pi_{\omega_2,1}(x)
 =\Pi_{\omega_2,1}(y)$
  \end{itemize}
 Let $C_{n,m,j}$ denote the image of $p_n$ under the mapping $Z_{m,j}$.
 \begin{itemize}
 \item[IH(9)] if $m < k \leq n$ , $j \in J_m$ and $i\in J_k$ then $C_{n,m,j}  \cap
 C_{n,k,i} = \emptyset$
 \end{itemize}
 By $[A_{m,j}, B_{m,j}]$ will be  denoted the interval whose endpoints
 are the two points $(f_m +
 H)(c^m_j)$ and $(f_m + H)(c^m_{j + 1})$. Observe that IH(5) implies that the definition
 of $[A_{m,j}, B_{m,j}]$ does not change at later stages of the induction.
 \begin{itemize}
 \item[IH(10)] the image of  $p_n$ under
  $ f_n\circ Z_{n,j} + \Phi_{n,j}\circ\Pi_{\omega_2,1}$ contains the interval
 $[A_{n,j}, B_{n,j}]$ for each $j \in J_n$
 \end{itemize}
 For $x \in [0,1]$ let $p_{n,m,j}^x$ be the join of all conditions $p_n\mid\sigma$ such that
 $\sigma : E_n \times  k_n \to 2$ is consistent with $p_n$ and $x$ belongs to the image of
 $(p_n\mid\sigma)^*$ under the mapping
 $f_n\circ Z_{m,j} + \Phi_{m,j}\circ \Pi_{\omega_2,1}$. The following is the key inductive requirement.
 \begin{itemize}
 \item[IH(11)] if $x \in [A_{m,j}, B_{m,j}]$ then
 $(p_{n+1,m,j}^x,k_{n+1}) <_{E_{n+1}}     (  p_{n,m,j}^x,k_n)$
 \end{itemize}

  Assuming that the induction can be completed, let $f
 =\lim_{n\rightarrow \infty} f_n$.  It will be shown that there is a
 condition $p_\omega\in \SS_{\omega_2}$ which forces that $T(f)$ belongs to
 the image of $N(f)$ under $f$. This contradiction will establish the
 theorem because IH(1) and IH(2) obviously guarantee that $f\in X$.

 Let $m$ be an integer such that there is some $j \in J_m$ such that $[
 c^m_j, c^m_{j+1}] \subseteq N(f)$ and $T(f)\in [A_{m,j}, B_{m,j}] $.
 The integers $m$ and $j$ must exist because the endpoints of $N(f)$
 belong to $D$ and so $N(f) = [c^m_i,c^m_k]$ for some $m$, $i$ and $k$.
 Furthermore, from IH(5) it follows that $[H+f(c^m_i), H+f(c^m_k)] =
 [H+f_m(c^m_i), H+f_m(c^m_k)] = \bigcup_{i\leq v < k} [A_{m,v},B_{m,v}
 ] $.  There must, therefore, be some $j$ between $i$ and
 $k-1$ which is suitable.

 It follows from IH(10) that the range of $f_m\circ Z_{m,j} +
 \Phi_{m,j}\circ\Pi_{\omega_2,1}$ contains $T(f)$ and so  $p^{T(f)}_{m,m,j}\neq \emptyset$.
 From IH(11) it follows that
 $$(p_{n+1,m,j}^{T(f)},k_{n+1}) <_{E_{n+1}}     (  p_{n,m,j}^{T(f)},k_n)$$
  for each $n \geq m$
 and so there is a condition in
 $p_\omega \in \SS_{\omega_2}$ such that  $p_\omega \leq
 p_{n,m,j}^{T(f)}$ for $n\geq m$.
 It follows that $T(f)$ belongs to the image of $p_\omega^* $ under the mapping
 $f\circ Z_{m,j} + \Phi_{m,j}\circ\Pi_{\omega_2,1}$.  Furthermore, because the
 diameters of the images of  $(p_{n,m,j}^{T(f)})^*$ under the mappings
 $f_n\circ Z_{m,j} + \Phi_{m,j}  $
 approach 0 as
 $n$ increases, it follows that $f\circ Z_{m,j} + \Phi_{m,j}\circ\Pi_{\omega_2,1}$ has constant value $T(f)$ on $p_\omega^*$.

 It follows that $p_\omega\forces{}{f\circ Z_{m,j} + \Phi_{m,j}\circ\Pi_{\omega_2,1}(G)
 =T(f)}$.  From IH(6) and the fact that $p_\omega \leq p_m$ it follows that
 $p_\omega\forces{}{H(z_{m,j}) = \Phi_{m,j}\circ\Pi_{\omega_2,1}(G)}$ and from IH(7) that
 $p_\omega\forces{}{z_{m,j} = Z_{m,j}(G)}$.  Therefore $p_\omega\forces{}{f(z_{m,j})
 + H(z_{m,j}) = T(f)}$ and this is a contradiction because $z_{m,j} \in
 [c^m_j, c^m_{j+1}] \subseteq N(f)$ by definition.

 To carry out the induction suppose that $f_n$, $\{\Phi_{n,j} \mid j\in
 J_n\}$ and $\{Z_{n,j}\mid j\in J_n\}$ as well as conditions $p_n \in
 \SS_{\omega_2}$ have all been defined
 for $n \leq K$. To begin, let $0 = c^{K+1}_0 < c^{K+1}_1 < c^{K+1}_2 <
 \ldots < c^{K+1}_{J_{K+1}} = 1$ be such that:
 \begin{itemize}
 \item $\{d_K\}\cup   {\cal C}_K \subseteq {\cal C}_{K+1}  =
 \{c^{K+1}_i \mid i\leq J_{K+1}\}$
 \item  the diameter of the image of $[c^{K+1}_i,c^{K+1}_{i+1}]$ under
 $f_K $ is less than $\epsilon_K\cdot 2^{-K-4}$
 \item $0 < \mid H(c^{K+1}_i) - H(c^{K+1}_{i+1})\mid < \epsilon_K\cdot 2^{-K-3}$
 \end{itemize}
 The first condition ensures that IH(4) is satisfied. The second is
 easily arranged using uniform continuity. The last condition can be
 satisfied by a further refinement using the Darboux property of $H$.

 Note that $\Pi_{\omega_2,1}(p^*_K)$, the image of $p_K^*$  under $\Pi_{\omega_2,1}$, is perfect and so,
 for each $i \in
 J_{K+1}$ it is possible to find $\Phi_{K+1,i} :\Pi_{\omega_2,1} (p_{K}^*)\rightarrow \RR$ such that
 \begin{itemize}
 \item $\Phi_{K+1,i}$ is a continuous mapping
 \item the image of $p_K^*$ under $\Phi_{K+1,i}\circ\Pi_{\omega_2,1}$ is  the interval whose endpoints
 are $H(c^{K+1}_i)$ and $H(c^{K+1}_{i+1})$
 \item  if $(m,j) \neq (K+1,i)$ then
 $\Phi_{K+1,i}\circ\Pi_{\omega_2,1}(x) \neq \Phi_{m,j}\circ\Pi_{\omega_2,1}(x)$  for every $x \in p_K^*$
 \item $\Phi_{K+1,i}$ is finite-to-one
 \end{itemize}
 Observe that the last point implies that $\Phi_{K+1,i}\circ\Pi_{\omega_2,1}(G)$ does not belong to the  ground
 model $V$.

 In any generic extension there must be a real between $c^{K+1}_i$ and
 $c^{K+1}_{i+1}$ at which $H$ takes on the value
 $\Phi_{K+1,i}\circ\Pi_{\omega_2,1}(G)$ because $H$ is assumed to be Darboux. Let
 $z_{K+1,i}$ be a name for such a real. It follows from the choice of
 $\Phi_{K+1,i}$ that $1\forces{}{z_{m,j}\neq z_{K+1,i}}$ for each
 $m\leq K+1$ and $j \in J_m$ such that $(K+1,i) \neq (m,j)$.

 Now find  $k$  and $p$ such that \begin{itemize}
 \item $(p,k_K) <_{E_{K+1}} (p_K, k_K)$
 \item $p$ is $(E_{K+1},k)$ determined
 \item  for each $\sigma : E_{K+1}\times k_{K}\to 2$ which is consistent with $p_K$
 and for each $m\leq K$, $j \in J_m$ and for each $x$ in the image of   $(p_K\mid\sigma)^*$
 under the mapping
 $f_K\circ Z_{m,j} + \Phi_{m,j}\circ\Pi_{\omega_2,1}  $     there is some $\sigma' :  E_{K+1}\times k\to 2$
 such that \begin{itemize}
 \item $\sigma'$ is consistent with $p$ \item
  $\sigma\subseteq \sigma'$ \item   the distance from $x$
 to the image of   $(p\mid\sigma')^*$
 under the mapping
 $f_K\circ Z_{m,j} + \Phi_{m,j}\circ\Pi_{\omega_2,1}  $     is less than $\epsilon_K\cdot 2^{-K-4}$\end{itemize}
 \item        the diameter of the image of   $(p\mid\sigma)^*$
 under the mapping
 $f_K\circ Z_{m,j} + \Phi_{m,j} \circ\Pi_{\omega_2,1} $     is less than $\epsilon_K\cdot 2^{-K-4}$ for each
 $\sigma  : E_{K+1}\times k\to 2$ which is consistent with $p$
 \end{itemize}
							     Now let  $k_{K+1}$ and $\bar{p}$ be such that
 \begin{itemize}
 \item $(\bar{p}, k_{K+1}) <_{E_{K+1}} (p,k)$
 \item $\bar{p}$ is $(E_{K+1},k_{K+1})$ determined
 \end{itemize}

 Because $V$ is closed under $H$ and $\Phi_{m,j}\circ\Pi_{\omega_2,1}  (G) \notin V$
 it follows that $z_{K+1,i}$ is a name for a real which does not belong
 to $V$. Lemma~\ref{name} can therefore be used $J_{K+1}$ times to find
 a condition $q$ such that $(q,{k_{K+1}}) \leq _{E_{K+1}}
 (\bar{p},{k_{K+1}})$ and for each $i\in J_{K+1}$ there is a function
 $Z_{K+1,i}: q^* \rightarrow [c^{K+1}_i,c^{K+1}_{i+1}]$ such that
 \begin{itemize}
 \item $q\forces{}{z_{K+1,i} = Z_{K+1,i}(G)}$
 \item $Z_{K+1,i}(x) \neq Z_{K+1,i}(y)$ unless $\Pi_{\omega_2,1}(x) = \Pi_{\omega_2,1}(y)$
 \end{itemize}
 Now observe that if $(m,j) \neq (K+1,i)$ then there can not be $\hat{q}
 \leq q$ such that $Z_{m,j}\restriction \hat{q}^*
 =Z_{K+1,i}\restriction \hat{q}^* $
 because it has already been remarked that
 $1\forces{}{z_{m,j}\neq z_{K+1,i}}$ for each
 $m\leq K+1$ and $j \in J_m$ such that $(K+1,i) \neq (m,j)$.
 It is therefore possible to  use Lemma~\ref{splitimage} repeatedly  to find
 a single  condition $p_{K+1}$  such that
 $(p_{K+1}, k_{K+1}) <_{E_{K+1}} ({q},{k_{K+1}})$ and
 the image of $p_{K+1}^*$ under $Z_{K+1,i}  $
 is disjoint from
 the image of $p_{K+1}^*$ under $Z_{m,j}  $
 if  $(m,j)\neq (K+1,i)$.      Observe that  $p_{K+1}$   is $(E_{K+1},k_{K+1})$
 determined because $\bar{p}$ is.  Hence IH(3) is satisfied.
 Now  define
 $C_{K+1,K+1,i}$ to be the range of $Z_{K+1,i}$.
 This, along with the induction hypothesis, will guarantee that  IH(6),
 IH(7), IH(8) and IH(9) are all  satisfied.

  For  integers
  $m \leq K+1$, $j \in J_m$ let $\{[u^0_{m,j,v},u^1_{m,j,v}] \mid v\in L_{m,j}\}$ be a
 partition of $[A_{m,j},B_{m,j}]$ into intervals of length $\epsilon_K\cdot 2^{-K- 2}$.
 Now, for each    $\sigma : E_{K+1}\times k_{K+1}\to 2$ and for each pair  of  integers
  $m \leq K+1$, $j \in J_m$  and for each
   $v\in L_{m,j}$  let
 $W[\sigma,m,j,v]$ be a
 perfect, nowhere dense subset of $$ (f_K + \Phi_{m,j}
 \circ(\Pi_{\omega_2,1}\restriction (p_{K+1}\mid\sigma))\circ
 Z_{m,j}^{-1})^{-1} [u^0_{m,j,v} - \frac{\epsilon_K}{2^{K + 3}  }   , u^1_{m,j,v} +
 \frac{\epsilon_K}{2^{K + 3}}] $$ if this is possible. By choosing smaller sets, if necessary,
  it may be assumed that the sets $W[\sigma,m,j,v]$ are pairwise disjoint and
 that $W[\sigma,m,j,v]\cap {\cal C}_{K+1} = \emptyset$.
 Then define $F_{\sigma,m,j,v} : W[\sigma,m,j,v]\rightarrow
 [u^0_{m,j,v},u^1_{m,j,v}]$ to be
 any continuous surjection and let $f_{\sigma ,m,j.v} = F_{\sigma ,m,j,v} -
 \Phi_{m,j}\circ \Pi_{\omega_2,1}\circ Z_{m,j}^{-1}$. Note that IH(8) implies that
 $\Pi_{\omega_2,1}\circ Z_{m,j}^{-1}$ is a function even though
 $Z_{m,j}$ is not one-to-one.

 Similarly, for each $i\leq
 J_{K+1}$ let $W_i$ be a perfect, nowhere dense subset of $[c^{K+1}_i,
 c^{K+1}_{i+1}]$ disjoint from each $ W[\sigma,m,j,v]$ and define $F_i : W_i
 \rightarrow [A_{K+1,i},B_{K+1,i}]$ to be a continuous surjection.
 Then let $f^i = F_i - \Phi_{K+1,i}
 \circ \Pi_{\omega_2,1}   \circ Z^{-1}_{K+1,i}$.  Notice
 that the domains of all the functions $f_{\sigma,m,j,v}$ and $f^i$ are
 pairwise disjoint.  Hence it is possible to find $f_{K+1} :I
 \rightarrow \RR$ extending each of these functions in such a way that
 $\rho( f_{K+1}, f_K)$  does not exceed  $$\max
 \{\mid f_{K+1}(y) - f_K(y)\mid  : y \in (\cup_iW^i) \cup(\cup_{\sigma,m,j,v}W[\sigma,m,j,v])\}$$
 and, moreover, because
 $W[\sigma,m,j,v] \cap {\cal C}_{K+1} = \emptyset$ and
 $W_i \cap {\cal C}_{K+1} = \emptyset$, it may also be arranged that
 $f_{K+1}(c) = f_K(c)$ if $c\in {\cal C}_K$. Therefore IH(5) is
 satisfied as well as IH(0).  Observe that IH(10) is satisfied because
  the choice of $F_j$ ensured
 that it maps $W_j$ onto $[A_{K+1,j}, B_{K+1,j}]$. Since
 $f_{K+1}\restriction W_j = F_j - \Phi_{K+1,j}  \circ \Pi_{\omega_2,1}\circ
 Z_{K+1,j}^{-1}$ it follows that $f_{K+1}\circ Z_{K+1,j} +
 \Phi_{K+1,j}\circ \Pi_{\omega_2,1}$ maps $p_{K+1} $ onto $
 [A_{K+1,j}, B_{K+1,j}]$.

 To see that IH(2) holds it suffices to consider only
 $$f_{K+1}\restriction (\cup_iW^i) \cup(\cup_{\sigma,m,j,v}W[\sigma,m,j,v])$$
 because $f_{K+1}$ was defined not to exceed this bound.
 Consider first  $y \in W[\sigma ,m,j,v]$. Then $\mid f_{K+1}(y)
 - f_K(y)\mid$ is equal to $$\mid f_{\sigma,m,j,v}(y) - f_K(y)\mid = \mid F_{\sigma,m,j,v}(y) -
 (\Phi_{m,j} \circ\Pi_{\omega_2,1} \circ Z_{m,j}^{-1}(y) + f_K(y))\mid$$ Next, it
 is possible to use the definition of $F_{\sigma,m,j,v}$ to conclude
  that $F_{\sigma,m,j,v}(y) \in
 [u^0_{m,j,v},u^1_{m,j,v}]$.   Using the definition of $W[\sigma,m,j,v]$ it is possible to conclude that
 $$(f_K + \Phi_{m,j} \circ\Pi_{\omega_2,1}  \circ
 Z_{m,j}^{-1})(y) \in [u^0_{m,j,v} - \frac{\epsilon_K}{ 2^{-K-3}},u^1_{m,j,v}+ \frac{\epsilon_K}{ 2^{-K-3}}]$$
  because $y \in W[\sigma,m,j,v]$.
 Consequently, $\mid f_{K+1}(y)
 - f_K(y)\mid   $ is no greater than the diameter of
 $$[u^0_{m,j,v} - \frac{\epsilon_K}{ 2^{-K-3}},u^1_{m,j,v}+ \frac{\epsilon_K}{ 2^{-K-3}}]$$  which is
 $\epsilon_K\cdot 2^{-K-3 } + \epsilon_K\cdot 2^{-K-2}   + \epsilon_K\cdot 2^{-K-3 } = \epsilon_K\cdot 2^{-K-1}$.

 On the other hand,
 if $y \in W_i$ then,
 as before, $$\mid f_{K+1}(y) - f_K(y)\mid =
 \mid F_{i}(y) -  (\Phi_{m,j} \circ\Pi_{\omega_2,1}  \circ Z_{m,j}^{-1}(y)  +
 f_K(y))\mid$$  Recall that  $\Phi_{K+1,i}$ is chosen to map onto $[H(c^{K+1}_i),
 H(c^{K+1}_{i+1})]$; moreover, because $y \in [c^{K+1}_i, c^{K+1}_{i+1}]$ it
 follows from the choice of ${\cal C}_{K+1}$ that $$f_K(y) \in
 [f_K(c^{K+1}_i) - \frac{\epsilon_K}{ 2^{-k-4}} , f_K(c^{K+1}_{i} + \frac{\epsilon_K}{ 2^{-k-4}}]$$
 and so
 $\Phi_{m,j} \circ\Pi_{\omega_2,1}  \circ Z_{m,j}^{-1}(y) + f_K(y)$ belongs to
 $[A_{K+1,i} - \frac{\epsilon_K}{ 2^{-K-4}}, B_{K+1,i}+\frac{\epsilon_K}{ 2^{-K-4}}]$.
   Furthermore, $F_{i}(y)$ belongs to $[A_{K+1,i}, B_{k+1,i}]$ by
 design.
 By the choice of ${\cal C}_{K+1}$  the
 diameter of $[A_{K+1,i}, B_{k+1,i}] $ is less than
 $\frac{\epsilon_K}{2^{-K-4}} +  \frac{\epsilon_K}{ 2^{-K-3}} $ and         so the diameter of
 $$[A_{K+1,i} - \frac{\epsilon_K}{ 2^{-K-4}}, B_{K+1,i}+\frac{\epsilon_K}{ 2^{-K-4}}]$$ is no greater  than
 $\frac{\epsilon_K}{ 2^{-K-1}}$ and
 so  it follows that $\mid f_{K+1}(y) - f_K(y)\mid < \epsilon_K\cdot 2^{-K-1}$.

 Now all of the induction hypotheses have been shown to be satisfied
 except for  IH(1) and IH(11).
 To verify  IH(11) suppose that
  $m \leq K$, $j\in J_m$ and $x\in [A_{m,j}, B_{m,j}]$. It follows that
 there is some $v\in L_{m,j}$ such that $x\in [u^0_{m,j,v}, u^1_{m,j,v}]$.
   Suppose also that
  $\sigma :E_{K+1}\times k_K\to 2$ is consistent
  with $p^x_{K,m,j}$. It follows that there is some   $\sigma' :E_{K+1}\times k\to
 2$ such that
 \begin{itemize}
 \item $\sigma'$ is consistent with $p$ \item
  $\sigma\subseteq \sigma'$
  \item   the distance from $x$
 to the image of   $(p\mid\sigma')^*$
 under the mapping
 $f_K\circ Z_{m,j} + \Phi_{m,j} \circ\Pi_{\omega_2,1}  $
 is less than $\epsilon_K\cdot 2^{-K-4}$\end{itemize}
  It suffices to show that  $\tau$ is consistent with  $p^x_{K+1,m,j}$
  for each  $\tau:E_{K+1}\times k_{K+1}\to 2$ such that
 $\sigma'\subseteq \tau$;
 the reason for this is that $k_{K+1}$ was chosen so that $(\bar{p},k_{K+1})
 \leq_{E_{K+1}} (p_K,k)$    and      $(p_{K+1},k_{K+1})
 \leq_{E_{K+1}} (\bar{p},k_{K+1})$.
  Recall that
  the diameter of the image of   $(p\mid\sigma')^*$
 under the mapping
 $f_K\circ Z_{m,j} + \Phi_{m,j}\circ\Pi_{\omega_2,1}   $     is less than $\epsilon_K\cdot 2^{-K-4}$ because
 $\sigma'  : E_{K+1}\times k\to 2$ is consistent with $p$.
  Because
  the distance from $x$ to
 to the image of   $(p\mid\sigma')^*$
 under the mapping
 $f_K\circ Z_{m,j} + \Phi_{m,j}\circ\Pi_{\omega_2,1}   $
 is less than $\epsilon_K\cdot 2^{-K-4}$
  it must be that  this image is
 contained in $[u^0_{m,j,v} - \epsilon_K\cdot 2^{-K-3}   , u^1_{m,j,v}
 - \epsilon_K\cdot 2^{-K-3}]$. Because $p_{K+1} < p$ it follows
 that the image of $(p_{K+1}\mid \tau)^*$ under the mapping $f_K\circ Z_{m,j} +
 \Phi_{m,j}\circ \Pi_{\omega_2,1}$ is contained in   $[u^0_{m,jv} - \epsilon_K\cdot 2^{-K-3}   , u^1_{m,jv} - \epsilon_K\cdot
 2^{-K-3}]$ and so $W[\tau,m,j,v]\neq \emptyset$.
 The choice of $F_{\tau,m,j,v}$ ensures that it maps    $W[\tau,m,j,v]$ onto
 $[u^0_{m,j,v}  , u^1_{m,j,v}]$ and therefore
 $f_{K+1} + \Phi_{m,j}\circ (\Pi_{\omega_2,1}\restriction (p_{K+1}\mid \tau)) \circ Z_{m,j}^{-1}$
 maps $W[\tau,m,j,v]$ onto   $[u^0_{m,j,v}  , u^1_{m,j,v}]$.
 Hence $\tau$ is consistent with $p^x_{K+1,m,j}$.

 Finally, choose $\epsilon_{K+1}$ so that the neighbourhood of $f_{K+1}$ of radius $\epsilon_{K+1}$
 is contained in $U_{K+1}$.
   \stopproof
  \begin{corol}
  If set theory is consistent \label{maincor}then it is consistent that for every Darboux function $F$ there is
		  a nowhere constant continuous function $f$ such that $F+f$ is also Darboux.
  \end{corol}
  \proof The model to use is the one for Theorem~\ref{main}. Given a Darboux $F$
 to obtain a nowhere constant continuous $f$ use the fact that the set of nowhere constant function
 is comeagre in ${\cal C}(I,\RR)$.
 \stopproof

   \section{Further Remarks} It should be observed that the
 function $f$ in Corollary~\ref{main} has very few nice properties
 other than continuity. It is natural to ask the following question.
 \begin{quest} Is \label{q1}there a Darboux function $H:I \rightarrow
 \RR$ such that $H+f$
 is not Darboux for every non-constant, differentiable function $f$?\end{quest}
 The answer to  Question~\ref{q1} for functions with continuous derivative is positive.
 The same question can be asked with
 absolutely continuous in the place of differentiable.  One should
 recall that differentiable functions satisfy the property $T_1$ of
 Banach \cite{Saks}.  \begin{defin} A function $F: \RR \rightarrow \RR$
 satisfies $T_1$ if and only if the set of all $x$ such that
 $f^{-1}\{x\}$ is infinite has measure zero.\end{defin} Banach showed
 that differentiable functions satisfy $T_1$. Question~\ref{q1} is of
 interest for differentiable functions because Corollary~\ref{ft1} shows that
 a strengthening of $T_1$
 yields a positive theorem.

 Another potentially interesting direction to pursue would be to ask whether the size
 of the set of
 continuous functions in Theorem~\ref{main} can be increased.
 \begin{quest} Is there a Darboux function  $F$ such that the set of continuous functions $f$
 such that $F+f$ is Darboux is comeagre?\end{quest}
  \begin{quest} Is there a Darboux function  $F$ such that the set of continuous functions $f$
 such that $F+f$ is Darboux has measure one?\end{quest}

 \nocite{br.ce.onsu}
 In \cite{be.na.unbd} the authors consider not only sums of a Darboux function and a
 continuous function but also products and other algebraic constructions. It is not
 difficult to check that everything that has been established in this paper  for sums
 also holds for products, but it is not clear that this must always be so.
 \begin{quest}
 If there is  a Darboux function $F$ such that $F+g$ is not Darboux for every nowhere
 constant function $g$ must it also be the case that
 there is  a Darboux function $F$ such that $F\cdot g$ is not Darboux for every nowhere
 constant function $g$? What about the opposite implication?                  \end{quest}

\makeatletter \renewcommand{\@biblabel}[1]{\hfill#1.}\makeatother
\renewcommand{\bysame}{\leavevmode\hbox to3em{\hrulefill}\,}

\end{document}